\documentclass[10pt]{amsart}
\usepackage[latin1]{inputenc}

\usepackage{amssymb}
\usepackage{color}

\newcommand{\N}{\mathbb{N}}
\newcommand{\Z}{\mathbb{Z}}
\newcommand{\Q}{\mathbb{Q}}
\newcommand{\R}{\mathbb{R}}
\newcommand{\C}{\mathbb{C}}

\newcommand{\Proj}{\mathbb{P}}

\newcommand{\om}{\omega}

\newcommand{\jj}{\mathcal{J}}
\newcommand{\J}{\mathcal{J}}

\newcommand{\Aa}{\mathcal{A}}
\newcommand{\Pp}{\mathcal{P}}
\newcommand{\E}{\mathcal{E}}

\newcommand{\Tt}{\mathcal{T}}

\DeclareMathOperator{\Gr}{Gr}
\DeclareMathOperator{\PD}{PD}
\DeclareMathOperator{\rk}{rk}

\DeclareMathOperator{\U}{U}
\DeclareMathOperator{\SO}{SO}

\DeclareMathOperator{\GL}{GL}
\DeclareMathOperator{\AGL}{AGL}

\DeclareMathOperator{\Lie}{Lie}
\DeclareMathOperator{\Ham}{Ham}

\DeclareMathOperator{\Symp}{Symp}
\DeclareMathOperator{\Diff}{Diff}

\DeclareMathOperator{\Aut}{Aut}

\newcommand{\CP}{\mathbb{CP}}

\newcommand{\PbP}{\CP^2\#\,\overline{\CP}\,\!^2}

\newcommand{\into}{\hookrightarrow}

\newcommand{\nblowup}[1]{\#\,#1\overline{\CP}\,\!^2}

\newcommand{\quotient}[2]{\raisebox{1pt}{\ensuremath{#1}}/\raisebox{-1pt}{\ensuremath{#2}}}

\def\eoe{\unskip\ \hglue0mm\hfill$\between$\smallskip\goodbreak}

\theoremstyle{plain}
\newtheorem{thm}{Theorem}[section]
\newtheorem*{thm*}{Theorem}
\newtheorem{prop}[thm]{Proposition}
\newtheorem*{prop*}{Proposition}
\newtheorem{lemma}[thm]{Lemma}
\newtheorem*{lemma*}{Lemma}
\newtheorem{cor}[thm]{Corollary}
\newtheorem*{cor*}{Corollary}

\theoremstyle{definition}

\newtheorem*{defn*}{Definition}
\newtheorem*{ackn}{Acknowledgments}

\theoremstyle{remark}
\newtheorem{remark}[thm]{Remark}
\newtheorem*{remark*}{Remark}
\newtheorem{example}[thm]{Example}
\newtheorem{example*}{Example}

\newtheorem*{facts*}{Facts}

\newcommand{\comment}[1]{}
\newcommand{\mute}[2]{}
\newcommand{\printlabel}[1]{}
\newcommand{\printversion}{}

\begin{document}

%%%%%%%%%%%%%%%%%%%%%%%%%%%%%%%%%%%%%%%%%%%%%%%%%%%%%%%%%%%%%%%%%%%%%%%%%%%%%%%%
% FRONT PAGE
%%%%%%%%%%%%%%%%%%%%%%%%%%%%%%%%%%%%%%%%%%%%%%%%%%%%%%%%%%%%%%%%%%%%%%%%%%%%%%%%

% TITLE

\title{Maximal compact tori in the Hamiltonian groups of $4$-dimensional symplectic manifolds}
\author{Martin Pinsonnault}
\address{Department of Mathematics, University of Toronto, Toronto, 
Canada, M5S 3G3.}
\curraddr{Fields Institute, Toronto, Canada, M5T 3J1.}
\email{mpinsonn@fields.utoronto.ca}
\thanks{Partially supported by a NSERC Postdoctoral fellowship.}
\printversion

% ABSTRACT

\begin{abstract}
We prove that the group of Hamiltonian automorphisms of a symplectic $4$-manifold contains only finitely many conjugacy classes of maximal compact tori with respect to the action of the full symplectomorphism group. We also extend to rational and ruled manifolds a result of Kedra which asserts that, if $M$ is a simply connected symplectic $4$-manifold with $b_{2}\geq 3$, and if $\widetilde{M}_{\delta}$ denotes a blow-up of $M$ of small enough capacity $\delta$, then the rational cohomology algebra of the Hamiltonian group $\Ham(\widetilde{M}_{\delta})$ is not finitely generated. Both results are based on the fact that in a symplectic $4$-manifold endowed with any tamed almost structure $J$, exceptional classes of minimal symplectic area are $J$-indecomposable. Some applications and examples are given.
\end{abstract}

\maketitle

% HEADINGS

\pagestyle{myheadings}
\markboth{MAXIMAL COMPACT TORI IN HAMILTONIAN GROUPS}{MARTIN PINSONNAULT}

%%%%%%%%%%%%%%%%%%%%%%%%%%%%%%%%%%%%%%%%%%%%%%%%%%%%%%%%%%%%%%%%%%%%%%%%%%%%%%%%
\section{Introduction}
%%%%%%%%%%%%%%%%%%%%%%%%%%%%%%%%%%%%%%%%%%%%%%%%%%%%%%%%%%%%%%%%%%%%%%%%%%%%%%%%

Let $(M,\om)$ be a closed symplectic manifold of dimension $2n$. Given a time-dependent function $H:M\times\R\to \R$, define the associate Hamiltonian vector field $X_{t}^{H}$ by the formula $dH_{t}=\om(X_{t}^{H},-)$. A Hamiltonian automorphism $\phi$ of $(M,\om)$ is a symplectomorphism of $M$ which can be represented as the time-$1$ map of the flow of some Hamiltonian vector field $X_{t}^{H}$. Denote by $\Ham(M,\om)$ the group of all Hamiltonian automorphisms equipped with the standard $C^{\infty}$-topology. This is an infinite dimensional  Fréchet Lie group whose Lie algebra consists of the space of smooth functions on $M$ with zero mean. 

The Hamiltonian group shares two remarkable properties with compact Lie groups: its Lie algebra comes with a natural bi-invariant nondegenerate inner product given by
$$(f,g) = \int_{M} fg\om^{n}$$
while the group itself admits a bi-invariant Finsler metric, called the Hofer metric, induced by the $L_{\infty}$ norm on its Lie algebra.
For this reason, it is generally believed that the Hamiltonian group $\Ham(M,\om)$ is tamer than the diffeomorphism group $\Diff(M)$ and constitutes an intermediate object between compact Lie groups and more general diffeomorphism groups, see~\cite{Bl},~\cite{MD:LecturesSymp}, and~\cite{Re}.

To develop a better understanding of this principle, it is natural to look at effective Hamiltonian group actions on symplectic manifolds. Since any compact Lie group $G$ contains a unique conjugacy class of maximal tori, and because the algebraic and homotopical properties of $G$ are, to a large extent, determined by the action of the Weyl group $W_{T}=N_{T}/T$ on a fixed maximal torus $T$, the first step in this direction is to understand maximal Hamiltonian actions by tori or, in other words, to classify conjugacy classes of maximal compact tori in $\Ham(M,\om)$. Because the symplectomorphism group $\Symp(M,\om)$ acts by conjugation on $\Ham(M,\om)$, this problem naturally breaks into two parts: one can first try to classify conjugacy classes of maximal tori in $\Ham(M,\om)$ with respect to the action of $\Symp(M,\om)$ and, afterward, with respect to the action of $\Ham(M,\om)$ itself.

The case of $2$-dimensional manifolds is classical: the only effective Hamiltonian torus action on a surface is, up to symplectic conjugation, the standard orthogonal $S^1$-action on $S^2\subset\R^{3}$. By Smale's theorem, $\Symp(S^{2})$ retracts onto $\SO(3)$, showing that $\Symp(S^{2})=\Ham(S^{2})$ and  that the standard circle action behaves homotopically as $S^{1}\subset \SO(3)$.

The case of $4$-manifolds is more interesting since the Hamiltonian group of a $4$-manifold can have distinct conjugacy classes of maximal tori. For example, the product $S^{2}\times S^{2}$ equipped with a product symplectic form $\mu\sigma\oplus\sigma$, $\mu\geq 1$, admits $\lceil \mu\rceil$ inequivalent Hamiltonian $T^{2}$-actions, where $\lceil \mu\rceil$ is the smallest integer larger than or equal to $\mu$, see~\cite{Ka-Hirzebruch}. Moreover, by performing a sequence of symplectic blow-ups of appropriate capacities on $S^{2}\times S^{2}$, we can also create Hamiltonian $S^{1}$-actions that cannot be extended to an action of $T^{2}$, that is, maximal circles subgroups in $\Ham(M,\om)$. Indeed, by blowing up rational manifolds, we can construct $4$-manifolds with any given number of inequivalent maximal $S^{1}$- or $T^{2}$-actions. However, our first result shows that we can create, up to symplectic equivalence, at most finitely many such actions, namely,
 
%: Finiteness
\printlabel{Finiteness}
\begin{thm}\label{Finiteness}
The Hamiltonian group of a symplectic $4$-ma\-ni\-fold contains at most finitely many conjugacy classes of maximal tori with respect to the action of the full symplectomorphism group.
\end{thm}
 
In contrast, the groups of volume preserving diffeomorphisms of several rational symplectic $4$-manifolds contain infinitely many conjugacy classes of maximal tori. So, from this point of view, Hamiltonian groups of $4$-manifolds do form intermediate objects between Lie groups and diffeomorphisms groups. 

The proof of Theorem~\ref{Finiteness} relies on Delzant's classification of toric actions, on Karshon's classification of Hamiltonian circle actions on $4$-manifolds, and on the analysis of exceptional curves in symplectic $4$-manifolds. 
More precisely, recall that a homology class $E\in H_{2}(M;\Z)$ is a \emph{symplectic exceptional class} if it can be represented by a symplectic sphere $S$ with $c_{1}(S)=1$ and $S\cdot S=-1$. In Sections 3, we reduce Theorem~\ref{Finiteness} to the following lemma:

%: LemmePrincipal
\printlabel{LemmePrincipal}
\begin{lemma}\label{LemmePrincipal}
Let $(M,\om)$ be a symplectic $4$-manifold not diffeomorphic to $\PbP$. Then, for any choice of $\om$-tame almost complex structure~$J$, all symplectic exceptional classes of minimal symplectic area are represented by an embedded $J$-holomorphic sphere. \end{lemma}

The problem of classisfying maximal tori in $\Ham(M,\om)$ up to \emph{Hamiltonian} conjugation is much harder and remains largely open. Nevertheless, Theorem~\ref{Finiteness} allows us to prove the following partial result:
\printlabel{HamConjugacyClasses}
\begin{thm}\label{HamConjugacyClasses}
Let $(M,\om)$ be a toric $4$-manifold. Then the set of Hamiltonian conjugacy classes of maximal $2$-tori in $\Ham(M,\om)$ is finite iff $\pi_{0}(\Symp(M,\om))$ is~finite. 
\end{thm}
Because two compact subgroups in $\Ham(M,\om)$ are Hamiltonian conjugate iff they define homotopic actions, this theorem indicates that the homotopy theoretic properties of maximal tori in $\Ham(M,\om)$ are rather subtle. Interestingly, Lemma~\ref{LemmePrincipal} turns out to be a key ingredient in the following generalization of a result of Kedra~\cite{Ke} on the rational homotopy type of symplectomorphism groups of $4$-manifolds:
\printlabel{KedraGeneralise}
\begin{thm}\label{KedraGeneralise}
Let $(M,\om)$ be a compact simply connected $4$-manifold with $b_{2}(M)\geq 3$. Let $\imath:B_{\delta}\into M$ be a symplectic embedding of the standard ball of capacity $\delta$ and denote by $(\widetilde{M}_{\delta},\widetilde{\om}_{\delta})$ the corresponding symplectic blow-up. Then, for $\delta$ sufficiently small, the rational cohomology of the Hamiltonian group of $(\widetilde{M}_{\delta},\widetilde{\om}_{\delta})$ is infinitely generated as an algebra.
\end{thm}

Recently McDuff~\cite{McDuff-ExtensionKedraResult} recast Kedra's result in a more general framework which allows for various generalizations, some of them valid in higher dimensions. It would be interesting to know if, similarly, a version of Theorem~\ref{Finiteness} holds in dimensions greater than $4$. Note that since there is no hope to classify Hamiltonian torus actions in every dimensions, any general proof of such a theorem would have to rely solely on intrinsic properties of symplectic and Hamiltonian groups.

In another direction, we do not know if there exist symplectic $4$-manifolds with infinitely many inequivalent symplecic, non-Hamiltonian actions. In view of the recent work of Pelayo~\cite{Pe}, this problem seems now accessible. It is interesting to note that there is no analog of Theorem~\ref{Finiteness} in the (full) contact category. Indeed, Lerman constructed~\cite{Le:ContactCuts} an example of an overtwisted contact $3$-manifold whose contactomorphism group contains infinitely many distinct conjugacy classes of maximal tori. To our knowledge, however, no examples with tight contact structures are known.

Finally, we note that Bloch, El Hadrami, Flaschka, and Ratiu drew very interesting analogies between compact Lie groups and Hamiltonian groups by studying, from an analytical point of view, centralizers of maximal tori in the Hamiltonian groups of two-dimensional annuli $A(r)\subset\R^{2}$ and toric manifolds. In particular, their results suggest\,\footnote{Note that the proof of their main theorem (Theorem~5.7) is flawed, so that the generalization of their results to toric manifolds is still conjectural.} that a version of the Schur-Horn-Kostant convexity theorem may hold for the Hamiltonian group of toric manifolds.

\begin{ackn}
The author would like to thank Yael Karshon and Dusa McDuff for their interest in this work and for many useful comments and suggestions, and the Fields Institute for providing a stimulating environment. This research was partially funded by NSERC Grant BP-301203-2004.

\end{ackn}

%%%%%%%%%%%%%%%%%%%%%%%%%%%%%%%%%%%%%%%%%%%%%%%%%%%%%%%%%%%%%%%%%%%%%%%%%%%%%%%%
\section{Exceptional curves in symplectic $4$-manifolds}
%%%%%%%%%%%%%%%%%%%%%%%%%%%%%%%%%%%%%%%%%%%%%%%%%%%%%%%%%%%%%%%%%%%%%%%%%%%%%%%%

In this section, we prove Lemma~\ref{LemmePrincipal}. We start by reviewing the symplectic blow-up construction. Then we recall some facts from Taubes-Seiberg-Witten theory of Gromov invariants and from the  classification of extremal rays on projective surfaces. Throughout this section, we assume the reader familiar with the theory of $J$-holomorphic curves in symplectic manifolds as explained in~\cite{AL}, \cite{Gr}, and~\cite{MS:J-HolomorphicCurves2}.

\subsection{Exceptional curves and symplectic blow-ups}

Let $(M,\om)$ be a symplectic $4$-manifold and let $\E=\E(M,\om)$ be the set of all symplectic exceptional homology classes $E\in H_{2}(M,\Z)$, that is, classes that can represented by an embedded symplectic sphere of self-intersection $-1$. By the adjunction formula, any such class must satisfy $c_{1}(E)=1$, where $c_{1}$ is the first Chern class of $TM$ with respect to any tame almost complex structure $J$. Moreover,

%: SetOfExceptionalClasses
\printlabel{SetOfExceptionalClasses}
\begin{lemma}\label{SetOfExceptionalClasses} Let $(M,\om)$ be a symplectic $4$-manifold.
\begin{enumerate}
\item The set $\E$ depends only on the deformation class of $\om$. 
\item Given an exceptional class $E\in\E$, the subset $\jj_{E}\subset\jj(M,\om)$ of tamed almost complex structures for which there exist an embedded holomorphic sphere in class $E$ is open, dense and path-connected. 
\item The subset $\jj_{\E}$ of $J$'s for which all classes in $\E$ are represented by embedded $J$-holomorphic spheres is of second category in $\jj$.
\item Given an arbitrary $J\in\jj$, each symplectic exceptional class is represented by either an embedded $J$-holomorphic sphere or by a $J$-holomorphic cusp-curve\,\footnote{A $J$-holomorphic cusp-curve is a connected union of possibly multiply-covered $J$-curves. Although the compactification of the moduli space is best described using the finer notions of \emph{stable maps} and \emph{stable curves}, the simpler notion of cusp-curves is sufficient for our purpose.} of the form
$$S=m_{1}S_{1}\cup\cdots\cup m_{n}S_{n}\text{,\hspace*{5mm}}n\geq 2$$
where each $m_{i}S_{i}$ stands for a simple spherical component occurring with multiplicity $m_{i}\geq 1$.
\item Any two symplectic exceptional classes $E,E'\in\E$ intersect non-negatively.
\end{enumerate}
\end{lemma}
\begin{proof} Statements (1) and (2) are proved in~\cite{MD:Structure}. Statement (3) follows from (2) and Baire's theorem. Statement (4) is a corollary of (3) together with Gromov's compactness theorem. The last statement follows from (3) and from positivity of intersections for $J$-holomorphic curves.
\end{proof}

The following lemma also follows from Gromov's compactness theorem.
%: Properness
\printlabel{Properness}
\begin{lemma}\label{Properness}
The period map $P_{\om}:\E\to \R_{+}$ which associates to a class $E$ its symplectic area $\om(E)$ is proper. Consequently, $P_{\om}$ has discrete image and reaches a minimal value $\varepsilon(M,\om)>0$. In particular, the set 
\[\E_{\min}=\{E\in\E ~|~\om(E)=\varepsilon(M,\om)\}\]
of smallest exceptional classes is finite and nonempty whenever $\E$ is nonempty. 
\end{lemma}
For now on, we will always denote by $\varepsilon(M,\om)$ the minimal area of an exceptional symplectic sphere in $(M,\om)$.

\printlabel{SoftProperness}
\begin{remark}\label{SoftProperness}
When $b_{2}^{+}(M)=1$, Lemma~\ref{Properness} follows from a purely algebraic argument. Indeed, the intersection form defines a metric of signature $(1,b_{2}-1)$ on $H_{2}(M;\R)$. It is then easy to show that given a homology class $A\in H_{2}(M;\R)$ such that $A\cdot A\geq 0$, and given any real numbers $a\leq b$ and $0<p\leq q$, the set
\[
N(A,a,b,p,q)=
\{ B\in H_{2}(M;\R)~:~-q \leq B\cdot B \leq -p\text{~~and~~}a\leq A\cdot B\leq b\}
\]
is compact. Lemma~\ref{Properness} follows by setting $A=\PD(\om)$ and $p=q=1$.

The compactness of the sets $N(A,a,b,p,q)$ has another interesting corollary: the group $\Aut_{[\om]}$ of automorphisms of $H_{2}(M;\R)$ which preserves the lattice $H_{2}(M;\Z)$ and which fixes the homology class $[\om]$ is always finite. Indeed, any automorphism $\phi\in\Aut_{[\om]}$ is determined by its action on
some orthogonal basis $\{A_{1},\ldots,A_{n}\}$ of $H_{2}(M;\Z)$. Now, given any  $1\leq i\leq k$, the set
\[
\Aa_{i}=\{B\in H_{2}(M;\Z)~|~B\cdot B=A_{i}\cdot A_{i}\text{~and~}\om(B)=\om(A_{i})\}
\]
is compact, hence finite. The finiteness of $\Aut_{[\om]}$ follows from the fact that any $\phi\in\Aut_{[\om]}$ acts by permutations on the sets $\Aa_{i}$.\eoe
\end{remark}

Recall that given a symplectic embedding $\imath:B^{4}_{\delta}\into (M,\om)$ of the standard ball of radius $r$ and capacity $\delta=\pi r^{2}$ in a $4$-manifold $(M,\om)$, one can perform a symplectic blow-up of size $\delta$ centered at $m=\imath(0)$ by removing the image of the ball and identifying a neighborhood of the boundary with a neighborhood of the zero section $\Sigma$ of the universal bundle over $\CP^{1}$. This operation produces a new symplectic manifold $(\widetilde{M}_{\imath},\tilde{\om}_{\imath})$ whose homology is naturally identified with $H_{*}(M)\oplus\Z[\Sigma]$, and in which the exceptional divisor $\Sigma$ is a symplectic exceptional sphere of symplectic area $\delta$. Since the construction depends on a connected family of parameters, the symplectic blow-up is well defined up to isotopy\,\footnote{Two symplectic forms $\om_{0}$ and $\om_{1}$ on the same manifold $M$ are \emph{deformation equivalent} if they belong to a one-parameter family of symplectic forms $\om_{t}$. They are \emph{isotopic} if all the forms $\om_{t}$ belong to the same cohomology class. By Moser's Stability Theorem, isotopic forms are symplectomorphic.} and hence up to symplectomorphisms. The inverse process is the symplectic blow-down and consists in replacing a neighborhood of a symplectic exceptional sphere $S$ by a ball of capacity $\om(S)$. Again, the resulting symplectic manifold is well defined up to isotopy. By analogy with the complex category, we say that two symplectic manifolds are rationally equivalent if one can be obtained from the other by a sequence of symplectic blow-ups and blow-downs.

It is often useful to interpret the symplectic blow-up in a slightly different way, namely, as a procedure which defines, out of any symplectic embedding $\imath:B^{4}_{\delta}\into (M,\om)$, a symplectic form on a fixed smooth blow-up $\widetilde{M}$ of $M$. As explained in~\cite{MP}, this approach allows comparison of blow-ups obtained from different embeddings:

\begin{prop}[\cite{MP} and \cite{MS:SymplecticTopology} Lemma 7.18]\hfill
\begin{enumerate}
\item Given two embeddings $\imath_{i}:B^{4}_{\delta_{i}}\to (M,\om)$ of arbitrary sizes $\delta_{1}$ and $\delta_{2}$, the corresponding blow-ups are deformation equivalent.
\item Given a family of embeddings $\imath_{t}:B^{4}_{\delta}\to (M,\om)$ of same capacity $\delta$, the blow-up manifolds $(\widetilde{M},\tilde{\om}_{t})$ are symplectomorphic. Moreover, if the embeddings are all centered at the same point $m\in M$, then the blow-up manifolds are isotopic.
\end{enumerate}
\end{prop}

No general statement about the homotopy type of embeddings spaces of symplectic balls of fixed capacity is known at the moment. However, in the special case of manifolds rationally equivalent to ruled $4$-manifolds, McDuff proved the following result:
\printlabel{UniquenessBlowups}
\begin{prop}[McDuff~\cite{MD:Isotopie}]\label{UniquenessBlowups} Suppose $(M,\om)$ is rationally equivalent to $\CP^{2}$ or to a ruled manifold $S^{2}\to N\to \Sigma_{g}$. Let $p$ be a point of $M$ and let $B^{4}_{\delta}\subset\R^{4}$ be the standard ball of capacity $\delta$. Then the space of centered embeddings $\imath:(B^{4}_{\delta},0)\into (M,p)$ is path-connected. Consequently, all symplectic blow-ups of $(M,\om)$ of same capacities are isotopic.
\end{prop}

\begin{remark}
The previous proposition does not imply that all cohomologous symplectic forms on a blow-up $\tilde{M}$ are isotopic. Only those symplectic forms obtained from isotopic forms on $M$ are known to be isotopic.\eoe
\end{remark}

Because symplectic exceptional spheres in the same homology class are symplectically isotopic, the blow-down process is easier to analyse:

\printlabel{UniquenessBlowdown}
\begin{prop}[McDuff~\cite{MD:Structure}]\label{UniquenessBlowdown}
\begin{enumerate}
\item The manifold obtained from a symplectic $4$-manifold $(M,\om)$ by blowing down a symplectic exceptional sphere $S$ is determined up to symplectic isotopy by the homology classe $[S]\in H_{2}(M;\Z)$. 
\item Every symplectic $4$-manifold $(M,\om)$ covers a minimal symplectic manifold which may be obtained from $M$ by blowing down a maximal collection of disjoint exceptional spheres. Moreover, the symplectomorphism class of the minimal model is unique unless $(M,\om)$ is the $k$-fold blow-up of $\CP^{2}$ or of a ruled manifold.
\end{enumerate}
\end{prop}

For manifolds rationaly equivalent to $\CP^{2}$ or to a ruled manifold the symplectic minimal model is often not unique. Worse, the blow-down of different sets of exceptional spheres may lead to the same minimal symplectic manifold. For instance, every rational manifold $M_{k}=\CP^{2}\nblowup{k}$ with $k\geq 9$ has infinitely many non-symplectomorphic minimal models. Fortunately, the ambiguity in the choice of a symplectic minimal model can always be reduced to a finite number of choices by considering blow-downs along spheres representing exeptional classes of minimal area, namely,

\printlabel{ChainOfMinimalBlowdowns}
\begin{cor}\label{ChainOfMinimalBlowdowns}
Given a symplectic manifold $(M,\om)$ rationaly equivalent to $\CP^{2}$ or to a ruled manifold, there exist finitely many maximal families of orthogonal exceptional classes $\{E_{k},\ldots ,E_{1}\}\subset H_{2}(M;\Z)$ such that the corresponding chain of blow-ups manifolds
\[(M,\om)=(M_{k},\om_{k})\stackrel{E_{k}}{\to}(M_{k-1},\om_{k-1})\stackrel{E_{k-1}}{\to} \cdots \stackrel{E_{1}}{\to}(M_{0},\om_{0})
\]
verify $E_{i}\in\E_{\min}(M_{i},\om_{i})$.
\end{cor}
\begin{proof}
Let $(M_{k},\om_{k})=(M,\om)$ and choose a symplectic exceptional class  of minimal area $E_{k}\in\E_{\min}(M_{k},\om_{k})$. By Lemma~\ref{Properness}, there are only finitely many such classes. Now, Proposition~\ref{SetOfExceptionalClasses} implies that $E_{k}$ is represented by an embedded symplectic sphere $S$ which can be symplectically blow-down. By Proposition~\ref{UniquenessBlowdown}, the symplectomorphism type of the resulting symplectic manifold only depends on $\om_{k}$ and on the choice of the class $E_{k}$. The results follows by induction.
\end{proof}

In the following, we will say that a blow-up (resp. a blow-down) is \emph{minimal} if the exceptional divisor (resp. the blow-down sphere) represents an exceptional class of minimal area. 

\begin{remark}
Although the choice of exceptional classes $\{E_{k},\ldots ,E_{1}\}\subset H_{2}(M;\Z)$ in Corollary~\ref{ChainOfMinimalBlowdowns} is usualy not unique, their areas $\om(E_{k})\geq\cdots\geq\om(E_{1})$ only depend on the initial manifold $(M,\om)$. Moreover, for a sufficiently generic symplectic form $\om$ on $M$, for example, when the periods of $\om$ are linearly independent over $\Q$, each set $\E_{\min}(M_{i},\om_{i})$ contains a unique class. In that case, the chain of minimal blow-ups is completely canonical.\eoe
\end{remark}

\subsection{Gromov invariants} 
Recall that for a 4-dimensional symplectic manifold $(M,\om)$ and a class $A\in H_{2}(X;\Z)$ such that $k(A)=\frac{1}{2}(A\cdot A + c_{1}(A))\geq 0$, the Gromov invariant of $A$ as defined by Taubes\footnote{The original definition given by Taubes must be modified slightly to handle multiply covered curves and to have an equivalence between Gromov invariants and Seiberg-Witten invariants. See~\cite{MD:LecturesGr} and~\cite{LL:EquivalenceGr-SW} for a complete discussion.} in~\cite{Ta} counts, for a generic almost complex structure $J$ tamed by $\om$, the algebraic number of embedded $J$-holomorphic curves in class $A$ passing through $k(A)$ generic points. In general, these curves may be disconnected and, in that case, components of negative self-intersection are exceptional spheres. An important corollary of Gromov's compactness theorem is that, for any choice of tamed $J$, a class $A$ with non zero Gromov invariant $\Gr_{\om}(A)$ is always represented by a union of $J$-holomorphic cusp-curves. In particular $\Gr_{\om}(A)\neq 0$ implies $\om(A)>0$.

The Gromov invariants only depend on the deformation class of $\om$, that is, $\Gr_{\om_{t}}(A)$ is constant along any $1$-parameter family $\om_{t}$ of symplectic forms. In fact, in dimension four, it follows from the work of Taubes that these invariants are smooth invariants of $M$. 
 
Now suppose that $(M,\om)$ is the $k$-fold blow-up of a ruled manifold 
$$S^{2}\to M_{0}\to\Sigma_{g}$$
Let $K_{M}$ stands for the canonical divisor of $(M,\om)$, that is, $K_{M}=-\PD(c_{1})$. Since $b_{2}^{+}=1$, the wall-crossing formula of \comment{Inserer les initiales} Li-Liu for Seiberg-Witten invariants (see~\cite{LL}, Corollary 1.4) implies that for every class $A\in H_{2}(M,\Z)$ verifying the condition $k(A)=\frac{1}{2}(A\cdot A + c_{1}(A))\geq 0$ we have
\begin{equation}\label{WallCrossing}
\Gr(A)\pm \Gr(K_{M}-A) = \pm (1+F\cdot A)^{g}
\end{equation}
where $F$ denotes the image in $H_{2}(M;\Z)$ of a fiber in $M_{0}$. 

\printlabel{NonZeroGr}
\begin{lemma}\label{NonZeroGr}
Suppose $(M,\om)$ can be obtained from a ruled manifold $(M_{0},\om_{0})$ by a sequence of symplectic blow-ups. Then any homology class $A\in H_{2}(M;\Z)$ verifying $A\cdot A\geq 0$, $\om(A)>0$, and $k(A)=\frac{1}{2}(A\cdot A + c_{1}(A))\geq 0$ has nonzero Gromov invariant.
\end{lemma}
\begin{proof}
Let $F_{0}\in H_{2}(M,\Z)$ be the image of a fiber of $(M_{0},\om_{0})$ under the natural inclusion $H_{2}(M_{0},\Z)\subset H_{2}(M,\Z)$ determined by the blow-up sequence defining $(M,\om)$. Note that we necessarily have $\om(F_{0})>0$.
Because $b_{2}^{+}(M)=1$, the set $\overline{\Pp}\in H_{2}(M;\Z)$ of integral homology classes with nonnegative self-intersection is a cone and $\overline{\Pp}-\{0\}$ has exactly 2 components. By the light-cone lemma, two classes in $\overline{\Pp}-\{0\}$ are in the same component if and only if they intersect nonnegatively. Since we have $\PD(\om)^{2}>0$, $F_{0}^{2}=0$, $\om(F_{0})>0$, and $\om(A)>0$, the classes $F_{0}$ and $A$ must lie in the component containing $\PD(\om)$. Consequently, $F_{0}\cdot C\geq 0$.

Because $k(A)\geq 0$, we can now apply the wall-crossing formula~(\ref{WallCrossing}) for Gromov invariants to deduce that $ \Gr(A)\pm \Gr(K-A) = \pm(1+F_{0}\cdot A) \neq 0$. In particular, $\Gr(A)$ and $\Gr(K-A)$ cannot be both equal to zero.  

We claim that $\Gr(K-A)= 0$. To see this, note that the invariance of Gromov invariants along deformations allows us to replace $\om$ by any convenient symplectic form in the same path-component. Note also that we can write the canonical class has $K=-\PD(c_{1})=-(B_{0}-E_{1}-\cdots -E_{k})$ where $B_{0}\in H_{2}(M_{0};\Z)$ is a class with positive symplectic area $\om_{0}(B)>0$ and where the $E_{i}$ are the exceptional classes defined by the blow-up sequence $M\to\cdots\to M_{0}$. By considering blow-ups of sufficiently small sizes, it is clear that we can find a deformation equivalent symplectic form $\om'$ for which the symplectic area of $K$ is strictly negative. Moreover, because $\om$ and $\om'$ are in the same component of $\overline{\Pp}-\{0\}$, any such form verifies $\om'(A)\geq 0$. Hence $\om'(K-A)<0$ and the Gromov invariant of $(K-A)$ must be zero. Consequently, $\Gr(A)=\pm(1+F_{0}\cdot A)^{g}\neq 0$.
\end{proof}

%: Mori theory of extremal curves
\subsection{Extremal curves on surfaces} Let us now explain some fundamental results from algebraic geometry. 

Recall that given an (integrable) complex structure $J$ on an algebraic variety $V$, the effective cone $NE^{J}\subset H_{2}(V,\R)$ is defined by setting
$$NE^{J}=\left\{A = \sum_{i}a_{i}[C_{i}]~~:~~C_{i}\text{~is a $J$-curve and $a_{i}\in\R_{+}$}\right\}$$
where $\R_{+}$ denotes the set of nonnegative real numbers. In particular, on a K\"ahler $4$-manifold $(M,\om,J)$, Lemma~\ref{SetOfExceptionalClasses} implies that $\E(M,\om)\subset\overline{NE}^{J}$. More generally, as was first noted by Ruan~\cite{Ru}, any class $A$ whose Gromov invariant is non zero must be contained in $\overline{NE}^{J}$.

We say that an element $A\in NE^{J}$ is \emph{extremal} if whenever $A$ can be written as a sum $A=B+C$ for $A, B\in NE^{J}$, then $B$ and $C$ are multiples of $A$. For smooth algebraic varieties, the structure of $NE^{J}(V)$ can be described in terms of extremal classes. To this end, given any cohomology class $\Gamma\in H^{1}(V;\Z)$, we define $(\Gamma)_{\leq 0}$ by setting
$(\Gamma)_{\leq0}=\left\{ A\in H_{2}(M;\R) | \Gamma(A)\leq 0 \right\}$.
Then, 
\printlabel{ConeTheorem}
\begin{thm}[Mori's Cone theorem~\cite{Mo}]\label{ConeTheorem}
Let $V$ be a smooth $n$-dimensional complex projective variety with first Chern class $c_{1}$. Then there exists a set of rational curves $l_{v}$ with $0<c_{1}(l_{i})\leq n+1$, called \emph{Mori extremal rational curves}, such that
\begin{enumerate}
\item $\overline{NE}^{J}(V)=\left(\sum\R_{+}[l_{i}]\right)+\overline{NE}^{J}(V)\cap(c_{1})_{\leq 0}$
\item Each $l_{i}$ is extremal.
\end{enumerate}
\end{thm}
 
The classification of extremal curves on complex projective surfaces is classical and can be traced back to the work of Castelnuovo and Enriques:
\printlabel{ClassificationExtremalCurves}
\begin{thm}[see Mori~\cite{Mo}]\label{ClassificationExtremalCurves} An irreducible curve C on an algebraic surface $V$ is an extremal rational curve if and only if one of the following conditions is satisfied:
\begin{enumerate}
\item $C$ is an exceptional curve,
\item $V$ has the structure of a ruled manifold such that $C$ is one of the fibers,
\item $V=\mathbb{P}^{2}$ and $C$ is a projective line.
\end{enumerate}
\end{thm}

\subsection{Exceptional classes of minimal symplectic area} The key fact in the proof of Lemma~\ref{LemmePrincipal} is the following simple observation:

% LemmeClassesPositives
\printlabel{LemmeClassesPositives}
\begin{lemma}\label{LemmeClassesPositives}
Let $(M,\om)$ be a non-minimal symplectic $4$-manifold rationally equivalent to a ruled manifold but not diffeomorphic to $\PbP$. Then any class $A\in H_{2}(M;\Z)$ verifying $c_{1}(A)\geq 1$, $A\cdot A\geq 0$, and $\om(A)>0$, must also satisfy $\om(A)\geq \varepsilon(M,\om)$.
\end{lemma}

\begin{proof}
Because $k(A)\geq 1$, $A\cdot A\geq 0$, and $\om>0$, we can apply Lemma~\ref{NonZeroGr} to deduce that the class $A$ has nonzero Gromov invariant. It follows that, for any choice of $\om$-tamed $J$, the class $A$ must be in the effective cone $NE^{J}$.

If the symplectic form $\om$ is integral, $(M,\om)$ can be realized as a complex projective variety $(M,J_{0},\om)\subset\CP^{N}$. Because $A\in \overline{NE}^{j_{0}}$, theorem~\ref{ConeTheorem} implies that $A$ can be written as a sum of the form
$$A=\sum a_{i}[l_{i}] + \text{~classes with non positive first Chern number}$$
But since $c_{1}(A)\geq 1$, at least one $J_{0}$-extremal class $[l_{i}]$ must enter this linear combination with strictly positive coefficient $a_{i}$. Because $(M,\om)$ is the blow-up of a ruled manifold, and because $(M,\om)$ is not itself ruled, theorem~\ref{ClassificationExtremalCurves} implies that this extremal class must be exceptional. It follows that $\om(A)>\varepsilon(M,\om)$. Because any rational class can be made integral by an appropriate rescaling, this finishes the proof in the case $[\om]$ is rational. 

To extend this argument to the case of a non rational symplectic form $\om$, it is sufficient to construct, for any $\delta>0$, a deformation $\om_{t}$ starting at $\om=\om_{0}$ such that
\begin{enumerate} 
\item $\om_{1}$ is rational,
\item $|\om_{1}-\om_{0}|_{C^{\infty}}<\delta$, and
\item the minimal area $\varepsilon_{t}=\varepsilon(M,\om_{t})$ of the symplectic exceptional spheres in $(M,\om_{t})$ varies continuously with $t$.
\end{enumerate}
To show that such a deformation always exists, choose $E_{\min}\in\E_{\min}(\om)$  and consider a maximal set of pairwise orthogonal exceptional classes $\{E_{\min},E_{1},\ldots,E_{k}\}$. By theorem~\ref{UniquenessBlowdown}, $(M,\om)$ is obtained by blowing up a minimal surface $(M_{0},\om_{0})$ along embedded balls of capacities 
\[
\om(E_{min}),\om(E_{1}),\ldots,\om(E_{k})
\] 
It follows that any exceptional class $E\in\E(M,\om)$ can be uniquely written as
\[
E=e_{B}B-e_{\min}E_{\min}-\sum_{i=0}^{k} e_{i}E_{i}
\] 
where $B\in H_{2}(M_{0},\Z)$ and $e_{B},e_{\min},e_{1},\ldots,e_{k}\in\Z$. Because distinct symplectic exceptional classes always intersect positively, the coefficients $e_{\min},e_{1},\ldots,e_{k}$ are all non-negative as long as $E$ is not equal to $E_{\min}$ or to one of the $E_{i}$. It follows that the symplectic area of such an exceptional class never decreases along a deformation of $\om$ in which the areas $\om_{t}(E_{\min})$ and $\om_{t}(E_{i})$ decrease. Therefore, one can achieve the desired deformation $\om_{t}$ by reducing the sizes of the blow-ups appropriately.
\end{proof}

\begin{proof}[Proof of Lemma~\ref{LemmePrincipal}] Let $(M,\om)$ be a non-minimal symplectic $4$-manifold not diffeomorphic to the ruled manifold $\PbP$. Suppose that for some almost complex structure $J\in\J(\om)$, a class $E_{\min}\in \E_{\min}(\om)$ has no embedded $J$-holomorphic representative. By Gromov's compactness, $E_{\min}$ must be represented by a $J$-holomorphic cusp-curve having at least 2 distinct components. Because $c_{1}(E_{\min})=1$, there must be a spherical simple component representing a homology class $C$ with $c_{1}(C)\geq 1$ and such that $\om(C)<\om(E_{\min})$. If $c_{1}(C)=1$, the adjunction inequality
\[
0\leq 1+\frac{C\cdot C-c_{1}(C)}{2}\leq\frac{C\cdot C+1}{2}
\]
implies that $C\cdot C\geq -1$, with equality if and only if $C$ is an exceptional class. The area of $E_{\min}$ being minimal among all symplectic exceptional classes, $C$ itself cannot be exceptional and hence $C\cdot C\geq 0$. Now, by McDuff~\cite{MD:Immersions} Corollary 4.2.1, we can always perturb $J$ and assume that the $J$-holomorphic representative of the class $C$ is immersed. Consequently, we can apply \cite{MD:ImmersedSpheres} Theorem 1.4 to conclude that $(M,\om)$ is the $k$-fold symplectic blow-up of a ruled manifold $S^{2}\to M_{0}\to\Sigma_{g}$. Hence, we can apply Lemma~\ref{LemmeClassesPositives} to conclude that $\om(C)\geq\varepsilon(M,\om)$. This contradiction shows that $E_{\min}$ can only be represented by embedded $J$-holomorphic curves.
\end{proof}

\begin{remark}
Consider the Hirzebruch surface of degree $2i+1\geq 3$, that is
\[
W_{2i+1}=\{ ([z_{0}:z_{1}],[w_{0}:w_{1}:w_{2}])\in\CP^{1}\times\CP^{2}~:~z_{0}^{2i+1}w_{1}=z_{1}^{2i+1}w_{0}\}
\]
The projection $\CP^1\times\CP^2\to\CP^1$ endows $W_{2i+1}$ with the structure of a K\"ahler $\CP^1$-bundle over $\CP^1$ which is diffeomorphic to $\PbP$ and in which the symplectic area of the exceptional divisor is larger than the area of the fiber class $F$. In these K\"ahler manifolds, the exceptional class $E$ is only represented by a holomorphic cusp-curve containing the zero section $s_{0}=\{ ([z_{0}:z_{1}],[0,0,1])\}$. Consequently, neither Lemma~\ref{LemmePrincipal} nor Lemma~\ref{LemmeClassesPositives} hold for $\PbP$.\eoe
\end{remark}

In many circumstances, the set $\E_{\min}(M,\om)$ can be described. For instance,
\printlabel{MinimalExceptionalClass}
\begin{lemma}\label{MinimalExceptionalClass}
Let $(M,\om)$ be a $4$-manifold and let $\widetilde{M}_{\delta}$ be a symplectic blow-up of capacity $\delta$. Let $\Sigma$ denote the exceptional divisor. Then, there exists a constant $\delta_{0}>0$ such that, for all $\delta< \delta_{0}$, the set $\E_{min}(\widetilde{M}_{\delta})$ contains only the class $[\Sigma]$.
\end{lemma}
\begin{proof}
It is sufficient to show that the symplectic area of an exceptional class $E\neq [\Sigma]$ never decreases as $\delta=\om(\Sigma)$ decreases. To see this, first write the class $E$ as a linear combination $E=A-s[\Sigma]$, where $A\in H_{2}(M,\Z)$ and $s\in\Z$. The symplectic area of $E$ is then equal to $\om(A)-s\delta$. But distinct symplectic exceptional classes always intersect positively, so that $0\leq E\cdot [\Sigma] = s$. The result follows readily. 
\end{proof}

\printlabel{OneThirdLemma}
\begin{lemma}\label{OneThirdLemma}
Consider $\CP^{2}$ with its standard Fubini-Study form normalized so that the symplectic area of a line is $1$. Let $(M,\om)$ be obtained from $\CP^{2}$ by $k$ blow-ups of sizes $1/3\geq \delta_{1}\geq\cdots\geq\delta_{k}$. Let $E_{1},\ldots,E_{k}$ be the homology classes of the corresponding exceptional divisors. Then a chain of minimal blow-ups is given by
\[
(M,\om)\stackrel{E_{k}}{\to}\cdots\stackrel{E_{2}}{\to}\PbP\stackrel{E_{1}}{\to} \CP^{2}
\]
\end{lemma}
\begin{proof}
Let $L$ be the class of a line. By positivity of intersections, every exceptional class $E$ in $\E(M,\om)\setminus \{E_{1},\ldots,E_{k}\} $ can be written
$$ E=e_{0}L-\sum e_{i}E_{i}$$
where the coefficients $e_{i}$ are nonnegative. Consequently,
\[
\om(E)=e_{0}-\sum e_{i}\delta_{i} \geq e_{0}-\sum_{i}\frac{e_{i}}{3} = \frac{c_{1}(E)}{3}=\frac{1}{3}
\]
By induction, this proves that $E_{i}\in\E_{\min}(M_{i},\om_{i})$. Moreover, because $\om(E_{1})<1/2$, the exceptional class in $(M_{1},\om_{1})\simeq \PbP$ is symplectically indecomposable, proving that we can extend the chain one step further to reach $\CP^{2}$. 
\end{proof}

%%%%%%%%%%%%%%%%%%%%%%%%%%%%%%%%%%%%%%%%%%%%%%%%%%%%%%%%%%%%%%%%%%%%%%%%%%%%%%%%
\section{Conjugacy classes of maximal tori}
%%%%%%%%%%%%%%%%%%%%%%%%%%%%%%%%%%%%%%%%%%%%%%%%%%%%%%%%%%%%%%%%%%%%%%%%%%%%%%%%

In this section, we prove the Finiteness theorem~\ref{Finiteness}. We first briefly recall standard results about Hamiltonian torus actions on $4$-manifolds. Again, we shall not give any details and we refer the interested reader to the papers~\cite{Au}, \cite{De}, \cite{Ka-Classification}, and the references therein. 

\subsection{Effective Hamiltonian torus actions on $4$-manifolds} 

Let suppose that a compact torus $T$ acts effectively on a symplectic $4$-manifold $(M,\om)$ in a Hamiltonian fashion. Since the orbits of a Hamiltonian $T$-action are isotropic and, on an open dense set, orbits are free, $T$ is either the circle $S^{1}$ or the two dimensional torus $T^{2}$. In particular, any effective Hamiltonian $T^{2}$-action on a $4$-manifold is maximal. Such an action is called \emph{toric}. We say that two torus actions are \emph{equivalent} if they differ only by the composition of a reparametrization $T^{k}\to T^{k}$ with a symplectomorphism $M\to M$. Note that the symplectomorphism group $\Symp(M,\om)$ acts by conjugation on the set of maximal tori in $\Ham(M,\om)$ and that two actions $\rho,\rho':T\to\Ham$ are equivalent if, and only if, their image subgroups  belong to the same orbit. For now on, we will always assume our torus actions both effective and Hamiltonian.

\subsection{Effective Hamiltonian $T^{2}$-actions} 

Recall that given a Hamiltonian $T^{2}$-action on $M$, the image of the moment map $\mu:M\to\mathfrak{t}^{*}\simeq\R^{2}$ is a convex polygon $\Delta$. The moment map preimage of an interior point is always a free orbit, while the preimage of an edge $e_{i}$ and of a vertex $v_{j}$ are, respectively, an invariant embedded symplectic sphere $S_{i}$ and a fixed point $p_{j}$. The moment polygon always satisfy the \emph{Delzant condition}: the two edges emanating from any vertex can be generated by a $\Z$-basis of the lattice $\Z^{2}$  (equivalently, the slopes of its edges are either rational or infinite, and every two consecutive edges have integral outward normal vectors $(x,y)$ and $(x',y')$ with $xy'-yx'=1$). Moreover, 

\begin{thm}[Delzant's Classification Theorem~\cite{De}]
The moment polygon determines the Hamiltonian space up to equivariant symplectomorphism. Conversely, given any rational polygon $\Delta$, one can construct a simply connected K\"ahler manifold $(M_\Delta,\omega_\Delta)$ with a Hamiltonian $T^{2}$-action which preserves the complex structure and whose moment map image is $\Delta$.
\end{thm}

Since the moment map is only defined up to an additive constant, the polygons $\Delta$ and $\Delta+C$ define the same toric manifold. Moreover, if two actions differ only by an automorphism of the torus, then their moment polygons are congruent modulo an element of the affine linear group $\AGL(2,\Z)=\GL(2,\Z)\ltimes\R^{2}$. Consequently, we have the following ``Delzant correspondence'':

%: DelzantCorrespondance
\printlabel{DelzantCorrespondance}
\begin{cor}\label{DelzantCorrespondance}
The set of Hamiltonian $T^{2}$-actions on $4$-manifolds, up to symplectic conjugation and reparametrization of $T^{2}$, is in bijection with the orbits of the action of $\AGL(2,\Z)$ on the set of Delzant polygons $\Delta\subset\R^{2}$.
\end{cor}

This correspondence reduces, in principle, the study of Hamiltonian toric $4$-manifolds to the combinatorics of Delzant polygons. In particular, the symplectic conjugacy classes of maximal $2$-dimensional tori in the group $\Ham(M,\om)$ can be enumerated provided we can identify which Delzant polygons, modulo AGL(2,Z)-congruence, give rise to manifolds symplectomorphic to $(M,\om)$.
 
The homological properties of the toric space $(M_{\Delta},\om_{\Delta})$ can be read directly on the moment polygon $\Delta$ once its edges has been cyclically ordered counterclockwise. To this end, let the \emph{size} of an edge be measured by its \emph{rational length}, which is characterized by being invariant under $\AGL(2,\Z)$-congruence and by being equal to the standard euclidean length along the coordinate axes. Then, 

\printlabel{HomologyDelzantPolygons}
\begin{prop}\label{HomologyDelzantPolygons} 
\begin{enumerate}
\item The number of edges of $\Delta$ is equal to $N=2+\rk H_2(M_\Delta)$ and the classes $[S_{1}],\ldots,[S_{N}]$ generate $H_{2}(M_{\Delta},\Z)$ over $\Z$. 
\item For $i\neq j$, the intersection number $[S_{i}]\cdot [S_{j}]$ is equal to the geometric intersection of the edges $e_{i}$ and $e_{j}$, while the self-intersection number of an invariant sphere is given by the formula 
\[ [S_{i}]\cdot [S_{i}]=\det\begin{bmatrix}n_{i+1}\\ n_{i-1}\end{bmatrix}\]
where $n_{i}$ denotes the primitive outward normal vector to the $i^{\text{th}}$ edge. By the adjunction formula, it follows that the sum $\sum_{i}[S_{i}]$ is equal to $\PD[c_{1}(TM_{\Delta})]$.
\item  The symplectic area of the sphere $[S_{i}]$ is $2\pi$ times the rational length of the edge $e_{i}$. In particular, the perimeter of $\Delta$ (measured in rational length) is equal to the pairing $\langle\omega_\Delta,c_{1}\rangle/2\pi$.
\item The symplectic volume of $M_\Delta$ is equal to the euclidean area of $\Delta$.
\end{enumerate}
\end{prop}

The Delzant condition together with the self-intersection formula of Proposition~\ref{HomologyDelzantPolygons}(2) imply that a Delzant polygon $\Delta$ with at least $5$ edges contains at least one edge, say $e_{i}$, whose moment map preimage is an invariant embedded sphere $S_{i}$ of self-intersection $-1$. The equivariant blow-down of such a sphere gives another toric manifold whose moment polygon $\Delta_{e_{i}}$ is obtained from $\Delta$ by gluing a triangle of area $ \om(S_{i})$ along the edge $e_{i}$. Therefore, any Delzant polygon with at least 5 edges can be equivariantly blow-down to a Delzant triangle or to a Delzant quadrilateral.  

The inverse process is the equivariant blow-up and consists in replacing an embedded and invariant symplectic ball $B$ of capacity $\delta$ centered at a fixed point $p\in M$ by an invariant exceptional sphere of area $\delta$. The moment map image of the resulting toric manifold $\tilde{M}_{\delta}$ is the polygon $\tilde{\Delta}_{\delta}$ obtained from $\Delta$ by removing a triangle of area $\delta/2\pi$ based at the vertex $\mu(p)$. 

Delzant's theorem implies
\printlabel{EquivariantBlowups} 
\begin{lemma}\label{EquivariantBlowups}
The equivariant symplectomorphism class of a $T^{2}$-equi\-va\-ri\-ant blow-down depends only on the homology class of the blow-down sphere. Similarly, the equivariant symplectomorphism class of a $T^{2}$-equi\-va\-ri\-ant blow-up depends only on its size $\delta$ and on the choice of the fixed point $p$.
Consequently, up to equivariant symplectomorphisms, there are only finitely many ways to perform an equivariant blow-up of given capacity $\delta$ on a toric manifold.
\end{lemma}

A Delzant triangle must have all its sides of equal rational length $\lambda$. Such a triangle corresponds to $\CP^{2}$, with its standard Fubini-Study form scaled such that $\om_{\lambda}(\CP^{1})=2\pi\lambda$, and equipped with the standard Hamiltonian action of $T^{2}$. 

Any Delzant quadrilaterals is $\AGL(2,\Z)$-equivalent to a \emph{Hirzebruch trapezoid}. Each Hirzebruch trapezoid is described by three positive numbers: its width $a\in\R_{+}$, its height $b\in\R_{+}$, and its slope $m\in\N$. These numbers must satisfy $a\geq b$ and $a>\frac{mb}{2}$. Conversely, for each triple $(a,b,m)$ verifying these inequalities corresponds a Hirzebruch trapezoid $\Delta(a,b,m)$. When $m$ is even, the manifold $X_{\Delta(a,b,m)}$ is symplectomorphic to $(S^{2}\times S^{2},\om_{ab})$, where $\om_{ab}$ is the product form for which the two $S^{2}$ factors have area $a$ and $b$ respectively. When $m$ is odd, $X_{\Delta(a,b,m)}$ is symplectomorphic to $(\PbP,\om_{ab})$ where $\om_{ab}$ is a form giving area $a$ to a projective line and area $b$ to the exceptional divisor. 

If we denote by $\lceil x\rceil$ the smallest integer greater or equal to $x$, the previous discussion can be summarized as follows:
\printlabel{ClassificationBasicActions}
\begin{prop}\label{ClassificationBasicActions}
\begin{enumerate}
\item A Delzant polygon is either a Delzant triangle or is obtained from a Hirzebruch trapezoid by a finite sequence of blow-ups. Consequently, a four-dimensional symplectic toric manifold is symplectomorphic to $\CP^2$ or is obtained from a Hirzebruch surface by a sequence of equivariant symplectic blow-ups. 
\item The standard Hamiltonian $T^{2}$-action on $\CP^{2}$ is unique   up to equivalence.
\item The symplectic ruled manifolds $(S^{2}\times S^{2},\om_{ab})$ and $(\PbP,\om_{ab})$ admit respectively $\lceil a/b\rceil$ and $\lceil a/b-1/2\rceil$ inequivalent Hamiltonian $T^{2}$-actions. 
\end{enumerate}
\end{prop}
 
\printlabel{SectionS1Actions} 
\subsection{Effective Hamiltonian $S^{1}$-actions}\label{SectionS1Actions}
\printlabel{SectionS1Actions} 

An effective Ha\-mil\-to\-nian $S^{1}$-actions on a symplectic $4$-manifold is characterized, up to equivariant symplectomorphisms, by a labelled planar graph $\Gamma$ which encodes the essential information about the orbit types of the action and the cohomology class of the symplectic form. More precisely, suppose the circle $S^{1}$ acts effectively on a symplectic $4$-manifold $(M,\om)$ in a Hamiltonian way. If the stabilizer $G_{p}$ of a point $p\in M$ is not trivial, then $G_{p}$ is either the circle $S^{1}$ or a cyclic subgroup $\Z_{k}$ of order $k\geq 2$. Consider the closure $O(G)$ of the set of points whose stabilizer is isomorphic to the subgroup $G\subset S^{1}$. From Morse theory applied to the moment map we know that:
\begin{enumerate}
\item Each component of the fixed point set $O(S^{1})$ is either a single point or a symplectic surface. The maximum and minimum of the moment map are each attained on exactly one component of the fixed point set, while fixed points on which the moment map is not extremal are isolated. In particular, the fixed point set contains at most two disjoint surfaces.
\item If $G=\Z_{k}$, $k\geq 2$, each connected component of $O(G)$ is a symplectic two-sphere on which the quotient circle $S^{1}/\Z_{k}$ acts with two fixed points.
\end{enumerate}

\smallskip
The graph associated to an effective Hamiltonian $S^{1}$-action on $(M,\om)$ is defined in the following way:
\begin{itemize}
\item To each component of the fixed point set corresponds a unique vertex.
\item Two vertices are connected by an edge if and only if the corresponding isolated fixed points are connected by a $\Z_{k}$-sphere.
\item Each vertex is labelled by the value of the moment map on the corresponding fixed point component. If an extremal vertex correspond to a symplectic surface $S$, two additional labels are attached: the genus of that surface, and its normalized symplectic area $\frac{1}{2\pi}\int_{S}\om$. 
\item Each edge is labelled by the isotropy weight $k$ of the corresponding $\Z_{k}$-sphere.
\end{itemize}

It follows from the definition of the graph that
\begin{itemize}
\item The fixed points are ordered by their moment map labels. If a $\Z_{k}$-sphere connects two fixed points $p$ and $q$ with moment map labels $\mu(p)>\mu(q)$, we call $p$ its \emph{north} pole and $q$ its \emph{south} pole.
\item For $k\geq 2$, a fixed point has isotropy weight $-k$ (resp. $k$) iff it is the north (resp. south) pole of a $\Z_{k}$-sphere.
\item Any vertex is reached by at most two edges. Because the action is effective, the two isotropy weights at a fixed point $p$ are relatively prime. In particular,  if an extremal vertex corresponds to a fixed surface, then it is not reached by any edge.
\item A non-free orbit is either a fixed point or belongs to a $\Z_{k}$-sphere.
\end{itemize}

Because the moment map of a circle action is only defined up to a constant, and because the reparametrization group of $S^{1}$ is isomorphic to $\Z_{2}$, the graph associated to a circle subgroup $S^{1}\subset\Ham(M,\om)$ is well-defined up to the action of $\AGL(1;\Z)\simeq \R\ltimes\Z_{2}$ on the moment map data.

The first part of Karshon's classification of Hamiltonian circle actions on $4$-manifolds is the following uniqueness theorem:

\printlabel{KarshonUniqueness}
\begin{thm}[Karshon~\cite{Ka-Classification}] The graph determines the Hamiltonian circle action up to equivariant symplectomorphism. 
\end{thm}

Obvious examples of circle actions are obtained by restricting an Hamiltonian $T^{2}$-action to a subcircle. The corresponding graph is obtained from the moment polygon $\Delta$ by projecting it on the corresponding subalgebra $\Lie(S^{1})\subset\Lie(T^{2})$, keeping track of the fixed points components, the moment map data, and the possible $\Z_{k}$-spheres. 

An other class of examples is provided by considering effective $S^{1}$-actions on ruled manifolds $\CP^{1}\to M_{g}\to \Sigma_{g}$, where $g\geq 1$ (in the case $g=0$, these are Hirzebruch surfaces). These actions are constructed by thinking of a ruled manifold $M_{g}$ as the projectivization of a bundle $L\oplus \C$, where $L\to\Sigma_{g}$ is a line bundle of degree $k$. The diffeomorphism type of the resulting manifold $M_{g}$ depends only on the parity of $k$: $M_{g}$ is diffeomorphic to the product $\Sigma_{g}\times S^{2}$ if $k$ is even, and to the non-trivial bundle $S^{2}\ltimes \Sigma_{g}$ if $k$ is odd. In all cases, the $S^{1}$-action induced by the standard action on the trivial factor $\C\to\Sigma_{g}$ defines a maximal Hamiltonian action on $\Proj(L\oplus \C)$ whose fixed point set consists in two holomorphic sections of self intersection $\pm k$.

As in the toric case, a Hamiltonian circle action can be equivariantly blow-up by some amount $\delta>0$ at a fixed point $p$, creating an invariant exceptional sphere $\Sigma$ of area $\delta$. Let $m\geq n$ be the isotropy weights at $p$. There are three cases to consider:
\begin{itemize}
\item If $m-n\geq 1$, and $mn\neq 0$, then $\Sigma$ is a $\Z_{m-n}$-sphere. 
\item If $m=n=\pm 1$, then $p$ is an extremal isolated fixed point and $\Sigma$ is a fixed extremal surface. 
\item If $p$ belongs to an extremal surface $S$ of area $\alpha$, the proper transform $\tilde{S}$ is a fixed surface of area $\alpha-\delta$ and $\Sigma$ is an invariant sphere on which $S^{1}$ acts freely. 
\end{itemize}
Because the graph determines the action up to equivariant symplectomorphisms, we have
\printlabel{S1EquivariantBlowups}
\begin{lemma}\label{S1EquivariantBlowups}
An $S^{1}$-equivariant blow-up depends only on its size $\delta>0$ and on the fixed point component where it is performed. Consequently, given a symplectic $4$-manifold endowed with a Hamiltonian circle action, there exist only finitely many ways to perform a equivariant blow-up of size $\delta$ on it. 
\end{lemma}
Moreover, we know from the graph if an equivariant blow-up of given size $\delta$ can be performed:
\printlabel{ConditionsS1BlowUps}
\begin{prop}[Karshon~\cite{Ka-Classification}]\label{ConditionsS1BlowUps}
An equivariant blow-up of size $\delta$ can be performed at a fixed point $p$ if, and only, if  the following conditions are satisfied:
\begin{enumerate}
\item If $p$ belongs to a fixed extremal surface $S$, then $\om(S)>\delta$.
\item If $p$ belongs to a $Z_{k}$-sphere $S$, then then $\om(S)>\delta$.
\item If $p$ is an interior fixed point, then $\min \mu<\mu(p)-\delta<\mu(p)+\delta<\max \mu$.
\item If $p$ is extremal and isolated then, for any other fixed point $q\neq p$, we have $|\mu(p)-\mu(q)|>\delta$. 
\end{enumerate} 
\end{prop}

The second part of Karshon's classification is the following characterization of action graphs:

\printlabel{KarshonExistence}
\begin{thm}[Karshon~\cite{Ka-Classification}]\label{KarshonExistence}
A planar graph $\Gamma$ is the graph of a Hamiltonian circle action on a symplectic $4$-manifold if, and only if, it can be obtained from the graph of an action defined on $\CP^{2}$ or on a ruled manifold $M_{g}$ by performing a sequence of equivariant blow-ups.
\end{thm}

The next proposition characterizes the graphs of \emph{non-maximal} circle actions on rational $4$-manifolds.

\printlabel{NonMaximalActions}
\begin{prop}[Karshon~\cite{Ka-Classification}]\label{NonMaximalActions}
An Hamiltonian $S^{1}$-action extends to a toric action if and only if each fixed surface has genus zero and each nonextremal level set for the moment map contains at most two non-free orbits. Moreover, any circle action whose fixed points are isolated extends to a toric action.
\end{prop}

Using this characterization of non-maximal actions, one can easily prove that
\printlabel{NoMaximalCircleActions}
\begin{lemma}[see Karshon~\cite{Ka-Hirzebruch}]\label{NoMaximalCircleActions}
On $\CP^{2}$, $(S^{2}\times S^{2},\om_{ab})$, and $(\PbP,\om_{ab})$, every Hamiltonian circle action extends to a toric action.
\end{lemma}

Finally, Karshon's Classification Theorem allows a complete understanding of maximal circle actions on minimal irrational ruled manifolds. Let normalize the symplectic form on a ruled manifold $M_{g}$ by fixing the symplectic area of its fiber $F$ to $1$ and by setting $\om(B)=\mu>0$, where $B$ is a section of self-intersection $0$ in the case of the trivial bundle, or a section of self-intersection $-1$ in the case of the non-trivial bundle. By the work of Lalonde-McDuff~\cite{LM}, it is known that, up to rescaling, every ruled $4$-manifold is symplectomorphic to one of these~$M_{g,\mu}$.

\printlabel{S1ActionsOnRuledManifolds}
\begin{cor}\label{S1ActionsOnRuledManifolds}
Let $M_{g,\mu}$ be a normalized ruled $4$-manifold with $g\geq 1$. Then $\Ham(M_{g,\mu})$ contains exactly $\lceil\mu\rceil$ inequivalent maximal $S^{1}$-actions, where $\lceil\mu\rceil$ is the smallest integer greater or equal to $\mu$.
\end{cor}
\begin{proof}
See the proof of Lemme 6.15 in Karshon~\cite{Ka-Classification}.
\end{proof}

\subsection{Proof of Theorem~\ref{Finiteness}} The classification theorems of Del\-zant and Karshon describe Hamiltonian torus actions up to equivariant symplectomorphisms. As such, they do not give the list of all maximal Hamiltonian torus actions on a given symplectic $4$-manifold. However,

%: CanonicalModel
\printlabel{CanonicalModel}
\begin{prop}\label{CanonicalModel}
Let $(M,\om)$ be the $k$-fold blow-up, $k\geq 1$, of a symplectic ruled $4$-manifold. Then there exists a symplectic ruled manifold $(M_{0},\om_{0})$ and an ordered set of $k$ disjoint exceptional classes $\{E_{1},\ldots,E_{k}\}\subset \E(M,\om)$ of areas $\delta_{1}\geq\ldots\geq\delta_{k}$ such that every effective Hamiltonian torus action on $(M,\om)$ is equivalent to an action obtained from an action on $(M_{0},\om_{0})$ through a sequence of $k$ equivariant blow-ups of capacities $\delta_{1},\ldots,\delta_{k}$.
\end{prop}

\begin{proof} By Corollary~\ref{ChainOfMinimalBlowdowns}, we know that there exists a chain of minimal symplectic blow-downs 
\[(M,\om)=(M_{k},\om_{k})\stackrel{E_{k}}{\to}(M_{k-1},\om_{k-1})\stackrel{E_{k-1}}{\to} \cdots \stackrel{E_{1}}{\to}(M_{0},\om_{0})
\]
ending with a \emph{ruled} manifold $(M_{0},\om_{0})$ and for which the symplectic areas $\delta_{k}=\om(E_{k})\leq\cdots\leq\delta_{1}=\om(E_{1})$ only depend on $(M,\om)$. Now, if $(M_{k},\om_{k})$ comes equipped with an effective Hamiltonian torus action $\rho_{k}$, choose a compatible almost-complex structure $J_{k}\in\jj(\om_{k})$ and consider the $\rho_{k}$-invariant almost complex structure $J_{k}'$ obtained by averaging the adapted metric $g_{k}(x,y)=\om_{k}(x,J_{k}y)$. By Lemma~\ref{LemmePrincipal}, the  class $E_{k}$ is represented by an invariant $J_{k}'$-holomorphic sphere which can be equivariantly blow-down. By Lemma~\ref{EquivariantBlowups}, this operation defines a Hamiltonian action $\rho_{k-1}$ on $(M_{k-1},\om_{k-1})$ whose conjugacy class in $\Symp(M_{k-1},\om_{k-1})$ depends only on $\rho_{k}$. The result follows by induction. 
\end{proof}

\printlabel{SmallBlowUpsOfCP2}
\begin{cor}\label{SmallBlowUpsOfCP2}
Consider $\CP^{2}$ with its standard Fubini-Study form normalized so that the symplectic area of a line is $1$. Let $(M,\om)$ be obtained from $\CP^{2}$ by $k$ blow-ups of sizes $1/3\geq \delta_{1}\geq\cdots\geq\delta_{k}$. Then all actions on $(M,\om)$ come from $\CP^{2}$ through equivariant blow-ups of sizes $\delta_{1}\geq\cdots\geq\delta_{k}$.
\end{cor}
\begin{proof}
This follows readily from Lemma~\ref{OneThirdLemma}.
\end{proof}

\printlabel{FinitenessToricCase}
\begin{cor}\label{FinitenessToricCase}
There are only finitely many inequivalent toric actions on any given symplectic $4$-manifold.
\end{cor}
\begin{proof}
Let $(M,\om)$ be a toric manifold. By Proposition~\ref{ClassificationBasicActions}~(1), $(M,\om)$ is symplectomorphic to a rational manifold. By Proposition~\ref{CanonicalModel}, every toric action on $(M,\om)$ comes from a single ruled manifold $(M_{0},\om_{0})$ through a sequence of $k$ equivariant blow-ups of fixed capacities $\delta_{1},\ldots,\delta_{k}$. By Proposition~\ref{ClassificationBasicActions}~(3), $(M_{0},\om_{0})$ admits only finitely many inequivalent toric actions, each of which corresponding to a Hirzebruch trapezoid. Then, Lemma~\ref{EquivariantBlowups} implies that there exist only finitely many ways to perform $k$ equivariant blow-ups of sizes $\delta_{1},\ldots,\delta_{k}$ on each of these trapezoids. Finally, Proposition~\ref{UniquenessBlowups} implies that all these $k$-fold blow-ups define toric manifolds symplectomorphic to $(M,\om)$.
\end{proof}

\printlabel{PreuveAlternative}
\begin{remark}\label{PreuveAlternative}
As mentioned in the Introduction, Corollary~\ref{FinitenessToricCase} can be proven by ``soft'' means, that is, without any references to $J$-holomorphic techniques. The main point is that any invariant exceptional sphere corresponds to a side of the polygon and, consequently, its symplectic area is bounded by $2\pi\langle\om,c_{1}\rangle$. Then, we can use Remark~\ref{SoftProperness} to prove that a toric $4$-manifold can be equivariantly blow-down to finitely many distinct Hirzebruch surfaces, see~\cite{KKP} for more details.

One drawback of this simpler proof is that it does not lead to an efficient algorithm to count the actual number of inequivalent toric actions on a given manifold. More importantly, this argument cannot be used in the case of $S^{1}$-actions since one can easily construct examples of rational $S^{1}$-manifolds for which the pairing $\langle\om,c_{1}\rangle$ is strictly negative.\eoe
\end{remark}

In the case of Hamiltonian $S^{1}$-actions, observe that in a chain of minimal blow-ups, a maximal $S^{1}$-action on $(M_{i},\om_{i})$ may well come from the equivariant blow-up of a $S^{1}$-action on $(M_{i-1},\om_{i-1})$ which is not maximal (in fact, Lemma~\ref{NoMaximalCircleActions} implies that this is always the case whenever $(M,\om)$ is rational). However,

\printlabel{BlowUpsNonMaximalActions}
\begin{lemma}\label{BlowUpsNonMaximalActions}
Let $(M,\om)$ be a symplectic $4$-manifold and denote by $\tilde{M}_{\delta}=(\tilde{M},\tilde{\om})$ its symplectic blow-up of size $\delta$. Then, up to equivalences, there exist at most finitely many non-maximal Hamiltonian $S^{1}$-actions on $(M,\om)$ which induce maximal $S^{1}$-actions on $\tilde{M}_{\delta}$.  
\end{lemma}
\begin{proof}
First, observe that if a $4$-manifold $(M,\om)$ admits a non-maximal Hamiltonian $S^{1}$-action, then it is rational. By Corollary~\ref{FinitenessToricCase}, there exist only finitely many inequivalent Hamiltonian $T^{2}$-actions on $(M,\om)$. Consequently, it is sufficient to show that at most finitely many maximal $S^{1}$-actions on $\tilde{M}_{\delta}$ comes from subcircles of a given toric action on $(M,\om)$.

Second, because a Hamiltonian $S^{1}$-action with only isolated fixed points always extends to a toric action, a maximal $S^{1}$-action must have at least one symplectic surface in its fixed point components. So, if we assume a maximal $S^{1}$-action $\tilde{\rho}$ is the equivariant blow-up of a non-maximal action $\rho$, then the description of $S^{1}$-equivariant blow-ups in section~\ref{SectionS1Actions} implies that either 
\begin{enumerate}
\item $\rho$ already has non-isolated fixed points, or 
\item a sphere of fixed points is created by blowing-up $\rho$ at an extremal isolated fixed point of equal weights $\pm 1$.
\end{enumerate}

In the first case, note that a subcircle $S^{1}\subset T^{2}$ has non-isolated fixed points if and only if it fixes pointwise one of the spheres corresponding to the sides of the moment polygon $\Delta$, that is, if and only if its Lie algebra $\Lie(S^{1})\subset\mathfrak{t}^{2}$ is perpendicular to one of the edges of $\Delta$. Clearly, only finitely many such subcircles exist, and Lemma~\ref{S1EquivariantBlowups} implies that each of them gives rise to at most finitely many maximal $S^{1}$-actions on $\tilde{M}_{\delta}$.

In the second case, note that the resulting $S^{1}$-action extends to a toric action since the blow-up operation itself can be made equivariant with respect to the whole torus $T^{2}$. Indeed, since we are blowing up at an extremal fixed point which is isolated, we can suppose that this fixed point is also extremal for the toric action and we only need to check that the capacity of the $S^{1}$-equivariant blow-up is the same as the capacity of the toric blow-up. But this follows from the fact that this fixed point has equal weights $\pm 1$.
\end{proof}

The proof of the finiteness theorem in the case of maximal Hamiltonian $S^{1}$-actions now proceeds, mutatis mutandis, as in the toric case: use Proposition~\ref{CanonicalModel} and Lemma~\ref{BlowUpsNonMaximalActions}, together with Karshon's classification theorem. This gives

%: FinitenessS1Case
\printlabel{FinitenessS1Case}
\begin{cor}\label{FinitenessS1Case}
A symplectic $4$-manifolds admits only finitely many inequivalent, maximal, Hamiltonian $S^{1}$-actions.
\end{cor}

This concludes the proof of Theorem~\ref{Finiteness}.\qed

\subsection{Counting conjugacy classes of maximal tori} Our proof of Theorem~\ref{Finiteness} gives a inductive procedure to count the number of distinct symplectic conjugacy classes of maximal tori in the Hamiltonian group of many symplectic $4$-manifolds $(M,\om)$:
\begin{enumerate}
\item Find a chain of minimal blow-ups $(M,\om)\stackrel{E_{k}}{\to}\cdots\stackrel{E_{1}}{\to} (M_{0},\om_{0})$.
\item List all maximal actions on $(M_{0},\om_{0})$ and perform all possible equivariant blow-ups of size $\delta_{1}$. Repeat for all other blow-ups of sizes $\delta_{2},\ldots, \delta_{k}$.
\item Eliminate redundancies in the resulting list of actions.
\end{enumerate}
Note that the first step can be easily completed if one starts with either a toric manifold whose Delzant polygon is known, or with the $k$-fold blow-up of a ruled manifold which is not rational. Indeed, in the former case, the minimal exceptional classes always correspond to sides of the polygon, while in the later case, the set of exceptional classes is always finite.

Applying this procedure to equal blow-ups of $\CP^{2}$, we can generalize a result of Karshon and~Kessler~\cite{KK}:
\printlabel{ManifoldsWithoutActions}
\begin{prop}\label{ManifoldsWithoutActions}
Consider $\CP^{2}$ with its standard Fubini-Study form normalized so that the symplectic area of a line is $1$. Let $(M,\om)$ be obtained from $\CP^{2}$ by $k$ blow-ups of size $\delta$. 
\begin{enumerate}
\item If $k\geq 4$, $(M,\om)$ does not admit any Hamiltonian $T^{2}$-action.
\item If $(k-1)\delta\geq 1$, $(M,\om)$ does not admit any Hamiltonian $S^{1}$-action.
\end{enumerate}
\end{prop}

To prove Proposition~\ref{ManifoldsWithoutActions}, we need the following lemma which follows easily from the description of $T^{2}$-equivariant blow-ups and from Proposition~\ref{ConditionsS1BlowUps}:
\printlabel{NecessaryConditions}
\begin{lemma}[Karshon-Kessler~\cite{KK}]\label{NecessaryConditions} Let normalize  $\CP^{2}$ so that $\om(\CP^{1})=1$. 
\begin{enumerate}
\item $\CP^{2}$ admits $k$ $T^{2}$-equivariant blow-ups of size $\delta>0$ iff $k\leq 3$ and $\delta<1/3$.
\item $\CP^{2}$ admits $k$ $S^{1}$-equivariant blow-ups of size $\delta>0$ iff $(k-1)\delta<1$.
\end{enumerate}
\end{lemma}

\begin{proof}[Proof of Proposition~\ref{ManifoldsWithoutActions}]
In view of Lemma~\ref{NecessaryConditions}, it is sufficient to show that, under the hypotheses of the Proposition, all actions on $(M,\om)$, if any, must come from $\CP^{2}$ thought equivariant blow-ups. For this to be true, we only need to find a chain of minimal blow-ups ending with $\CP^{2}$. There are several cases to consider depending on $k$. If $k\geq 9$, the volume condition $0<[\om]\cdot[\om]=1-k\delta^{2}$ implies that $\delta\leq 1/3$ and the conclusion follows from Lemma~\ref{SmallBlowUpsOfCP2}. For $1\leq k\leq 8$, the set $\E(M,\om)$ is finite and can be explicitely described in terms of the natural basis $\{L,E_{1},\ldots,E_{k}\}$ defined by the blow-ups, see~\cite{De}, For $4\leq k\leq 8$, the conclusion follows from a case-by-case study of exceptional classes. For $2\leq k\leq 3$, the conditions of statements (1) and (2) are empty because of Gromov's packing inequality (see~\cite{MP}), while there is nothing to prove for $k=1$. 
\end{proof}

\subsection{Proof of Theorem~\ref{HamConjugacyClasses}}
We now consider Hamiltonian conjugacy classes of $2$-tori in $\Ham(M,\om)$.
\begin{prop}
Consider an effective toric action $\rho:T^{n}\into \Symp(M^{2n},\om)$ with moment map $\mu:M\to\mathfrak{t}\simeq\R^{n}$. Denote by $T$ the corresponding torus in $\Symp(M,\om)$, by $C_{T}$ its centralizer, by $N_{T}$ its normalizer, and by $W_{T}=N_{T} / C_{T}$ its Weyl group.
\begin{enumerate}
\item  The normalizer $C_{T}$ is equal to the group of all symplectomorphisms $\phi$ that preserve the moment map, that is, such that $\mu\circ\phi=\mu$. Moreover, $C_{T}$ is connected and $C_{T}\subset\Ham(M,\om)$. 
\item The normalizer $N_{T}$ is equal to the group of all symplectomorphisms $\psi$ such that $\mu\circ\psi=\Lambda\circ\mu$ for some $\Lambda\in\AGL(n;\Z)$. In particular, the Weyl group $W_{T}$ is finite.
\end{enumerate} 
\end{prop}
\begin{proof} 
(1) Let $\phi$ be an element of $C_{T}$ and let $t$ be an element of $T$ generated by an Hamiltonian function $H:M\to\R$. Then, $\phi t \phi^{-1}$ is generated by $H\circ \phi^{-1}$. It follows that the moment map of the conjugated action $\phi\rho\phi^{-1}$ is $\mu\circ\phi^{-1}$. Since $\phi\rho\phi^{-1}=\rho$, there exists a constant $C\in\R^{n}$ such that $\mu\circ\phi^{-1}=\mu+C$. But because $M$ is compact, $C=0$, and $\phi$ preserves the moment map.

Because the action is toric, the fibers of the moment map are orbits of the action. It follows from a theorem of Schwartz~\cite{Sch} that $\phi$ can be written in terms of the action as $\phi(m)=\overline{f}(m)\cdot m$, where $\overline{f}:M\to T^{n}$ is a smooth function which is constant on each fiber. By Corollary 4.2 of~\cite{KL}, it follows that $\overline{f}$ is the pullback by $\mu$ of a smooth function $\Delta\to T^{n}$, and because $\Delta$ is contractible, there exists a smooth function $f:\Delta\to\mathfrak{t}$ such that $\overline{f}=\exp [f\circ\mu]$ (here, $\exp:\mathfrak{t}\to T^{2}$ denotes the exponential map). Now, the Hamiltonian flow of the real function $F$ defined by the pairing $F(m)=\langle \mu(m),f\circ \mu(m) \rangle$ is
\[
\phi_{s}(m)=\exp[ s (f\circ\mu(m))]\cdot m
\]
Since $\phi_{s}\in C_{T}$ and $\phi=\phi_{1}$, this shows that $C_{T}$ is connected and that $C_{T}\subset\Ham$.

\medskip
(2) Every element of the normalizer $N_{T}$ acts on $T$ by automorphisms. Consequently, given a parametrization of the torus, there is an exact sequence
$$1\to C_{T}\to N_{T}\to W_{T}\to 1 $$
in which $W_{T}\subset \GL(n;\Z)$ is discrete. Given $\psi\in N_{T}$ this implies that the moment map of $\psi\rho \psi^{-1}$ equals $\Lambda\circ \mu$ for some $\Lambda\in\AGL(n;\Z)$. On the other hand, the moment map of $\psi\rho\psi^{-1}$ is also equal to $\mu\circ \psi$, showing that 
\begin{equation}\label{EqnNormalizer}
\mu\circ\psi=\Lambda\circ\mu
\end{equation}
The finiteness of $W_{T}$ follows from the fact that, up to translations, there exist only finitely many $\Lambda$ verifying~(\ref{EqnNormalizer}). 
\end{proof}

Given a rational symplectic $4$-manifold $(M,\om)$, let $\Tt$ denotes the set of all $2$-tori in $\Ham(M,\om)$. Then, Theorem~\ref{HamConjugacyClasses} states that for any toric $4$-manifold $(M,\om)$, the set $\Tt/\Ham$ is finite iff $\pi_{0}(\Symp)$ is finite.

\begin{proof}[Proof of Theorem~\ref{HamConjugacyClasses}]
Suppose the set $\Tt/\Ham$ contains $n<\infty$ conjugacy classes. Because toric $4$-manifolds are simply connected, $\pi_{0}(\Symp)=\Symp/\Ham$. Consequently, the group $\pi_{0}(\Symp)$ acts on the set $\Tt/\Ham$ and there is an induced homomorphism $c:\pi_{0}(\Symp)\to\mathbb{S}_{n}$. Let $\pi_{T}$ be the isotropy subgroup of $[T]\in\Tt/\Ham$, that is
\begin{eqnarray*}
\pi_{T}&=&\{[g]\in\pi_{0}(\Symp)~|~\exists\, h\in\Ham \text{~such that~} hg\in N_{T}\}\\
&=&\quotient{N_{T}}{N_{T}\cap\Ham}
\end{eqnarray*}
Because $C_{T}\subset\Ham$, there is a quotient map $W_{T}\to\quotient{N_{T}}{N_{T}\cap\Ham}$. Consequently, $\pi_{T}$ is finite. The finiteness of $\pi_{0}(\Symp)$ now follows from the exact sequence $1\to\ker(c)=\bigcap_{T} \pi_{T}\to\pi_{0}(\Symp)\to\mathbb{S}_{n}$.

Conversely, if $\Tt/\Ham$ is infinite, the map $\Tt/\Ham\to\Tt/\Symp$ has at least one infinite fiber consisting of distincts classes $[T_{i}]$, $i\in\N$, for which there exist $g_{i}\in\Symp$ such that $g_{i}T_{i}g_{i}^{-1}=T_{0}$. Assume that for two indices $i\neq j$ we have $[g_{i}]=[g_{j}]$ in $\pi_{0}(\Symp)$. Then there exists $h\in\Ham$ such that $hg_{i}=g_{j} $ and $(g_{i}^{-1}h^{-1}g_{i})T_{i}(g_{i}^{-1}h^{-1}g_{i})^{-1}=T_{j}$. Since $(g_{i}^{-1}h^{-1}g_{i})\in\Ham$ this is impossible, showing that the isotopy classes $[g_{i}]\in\pi_{0}(\Symp)$ are all distinct.
\end{proof}

Let $\Symp_{h}(M,\om)$ denote the subgroup of symplectomorphisms that act as the identity on $H_{*}(M;\Z)$. Because for symplectic $4$-manifolds with Hamiltonian circle actions the group $\Aut_{[\om]}$ is finite (see the second part of Remark~\ref{SoftProperness}), the group $\pi_{0}(\Symp(M,\om))$ is finite iff $\pi_{0}(\Symp_{h}(M,\om))$ is finite. For $(M,\om)$ diffeomorphic to $\CP^{2}$, $S^{2}\times S^{2}$, $\PbP$, or $\CP^{2}\nblowup{2}$, we know from~\cite{AM,Gr,Pi} that $\Ham(M,\om)=\Symp_{h}(M,\om)$. Because there are no maximal $S^{1}$-actions on these manifolds, this gives
\begin{cor}\label{CorHamClasses}
There are at most finitely many Hamiltonian conjugacy classes of maximal tori in the Hamiltonian groups of a symplectic manifold diffeomorphic to $\CP^{2}$, $S^{2}\times S^{2}$, $\PbP$, or $\CP^{2}\nblowup{2}$.
\end{cor}

\begin{remark}
It is likely that Corollary~\ref{CorHamClasses} also holds for $\CP^{2}\nblowup{3}$ and $\CP^{2}\nblowup{4}$.
\end{remark}

\begin{remark}
Presently we do not know a single toric $4$-manifold $(M,\om)$ for which $\pi_{0}(\Symp)$ is infinite, although we believe that such manifolds should be quite common. The reason is that the only examples of rational $4$-manifolds for which $\pi_{0}(\Symp)$ is known to be infinite consist of \emph{monotone} $\CP^{2}\nblowup{k}$ with $k\geq 5$ and, by Proposition~\ref{ManifoldsWithoutActions}, these manifolds do not admit any hamiltonian circle action. These results are due to Seidel~\cite{Se} who found a necessary condition for iterated symplectic Denh twists along Lagragian spheres to be symplectically isotopic to the identity. This condition involves the algebraic structure of the quantum ring $QH_{*}(M;\Lambda_{\om}^{\text{univ}})$ with coefficients in the universal Novikov ring $\Lambda_{\om}^{\text{univ}}$, and it gives an effective criteria only when $(M,\om)$ is monotone. Although it is still not clear whether this limitation can be removed, it would be interesting to study some simple cases, such as the $5$-fold blow-up of $\CP^{2}$ with sizes $\{1/3,1/3,1/3,\epsilon,\epsilon\}$, where $\epsilon<1/3$.\eoe
\end{remark}

%%%%%%%%%%%%%%%%%%%%%%%%%%%%%%%%%%%%%%%%%%%%%%%%%%%%%%%%%%%%%%%%%%%%%%%%%%%%
\section{Rational cohomology of symplectomorphism groups}
%%%%%%%%%%%%%%%%%%%%%%%%%%%%%%%%%%%%%%%%%%%%%%%%%%%%%%%%%%%%%%%%%%%%%%%%%%%%

In \cite{Ke}, Kedra proved Theorem~\ref{KedraGeneralise} in the case $(M,\om)$ is not the $k$-fold blow-up of a rational surface. His main result can be restated in the following way:

\printlabel{Kedra}
\begin{thm}[Kedra~\cite{Ke}]\label{Kedra}
Let $(M,\om)$ be a compact simply connected $4$-manifold with $b_{2}(M)\geq 3$. Let $\imath:B_{\delta}\into M$ be a symplectic embedding of the standard ball of capacity $\delta$ and denote by $(\widetilde{M}_{\delta},\widetilde{\om}_{\delta})$ the corresponding symplectic blow-up. Suppose that for small enough $\delta$ the exceptional divisor $\Sigma\subset (\widetilde{M}_{\delta},\widetilde{\om}_{\delta})$ is always represented by a $\tilde{J}$-holomorphic representative, for any choice of tamed almost-complex structure $\tilde{J}$ on $(\widetilde{M}_{\delta},\widetilde{\om}_{\delta})$. Then, for $\delta$ sufficiently small, the rational cohomology of the Hamiltonian group of $(\widetilde{M}_{\delta},\widetilde{\om}_{\delta})$ is infinitely generated as an algebra.
\end{thm}

The proof, which we outline below, combines known facts from homotopy theory together with elementary properties of symplectic balls in $4$-manifold, see~\cite{Ke} and~\cite{LP-Duke} for more details. 

The first step is to apply rational homotopy techniques and Gottlieb theory of evaluation subgroups to prove the following topological results: 
\begin{itemize}
\item If a topological group acts transitively on a simply connected \emph{rationally hyperbolic}\,\footnote{A topological space $M$ is rationally hyperbolic if the dimensions of its rational homotopy groups $\pi_{k}(M)\otimes\Q$ grow exponentially.} space $M$, then the rational cohomology of the isotropy subgroup of a point is infinitely generated as an algebra.
\item Every simply connected $4$-manifold $M$ for which $\dim H_{2}(M;\Q)\geq 3$ is rationally hyperbolic.
\end{itemize}
In particular, given any symplectic, simply connected $4$-manifold $(M,\om)$ with $b_{2}\geq 3$, the rational cohomology of the stabilizer $\Symp(M, p;\om )$ is infinitely generated as an algebra. 
 
Now consider the subgroup $\Symp^{\U(2)}(M,B;\om)$ of symplectomorphisms that restrict to a unitary transformation on a given standard symplectic ball $B_{\delta}$ of capacity $\delta$ centered at $p$. The second step is to show that, if $\delta$ is small enough, the inclusion induces a map $\pi_{*}(\Symp^{\U(2)}(M,B;\om))\otimes\Q\to\pi_{*}(\Symp(M, p;\om ))\otimes\Q$ whose image has infinite rank.

Finally, because every element of $\Symp^{\U(2)}(M,B;\om)$ lifts to the blow up $(\tilde{M_{\delta}},\tilde{\om}_{\delta}) $, the group $\Symp^{\U(2)}(M,B;\om)$ maps into the subgroup $\Symp(\tilde{M_{\delta}},\Sigma;\tilde{\om}_{\delta})$ of symplectomorphisms that fix the exceptional divisor $\Sigma$. The last step then consists in showing that if the exceptional class $[\Sigma]$ is symplectically indecomposable, each map in the sequence
$$\Symp^{\U(2)}(M,B;\om)\simeq \Symp(\tilde{M_{\delta}},\Sigma;\tilde{\om}_{\delta})\simeq \Symp(\tilde{M_{\delta}};\tilde{\om}_{\delta})$$
is a homotopy equivalences, see~\cite{LP-Duke}. This shows that the rank of the module $\pi_{*}(\Symp(\tilde{M_{\delta}};\tilde{\om}_{\delta}))\otimes \Q$ is infinite and that the rational cohomology of $\Symp(\widetilde{M}_{\delta},\widetilde{\om}_{\delta})$ is infinitely generated as an algebra.

\begin{proof}[Proof of Theorem~\ref{KedraGeneralise}] Assuming Theorem~\ref{Kedra}, the proof of Theorem~\ref{KedraGeneralise} reduces to showing that on any compact simply connected $4$-manifold with $b_{2}(M)\geq 3$, and for small enough $\delta$, the exceptional divisor $\Sigma\subset (\widetilde{M}_{\delta},\widetilde{\om}_{\delta})$ is always represented by a $\tilde{J}$-holomorphic representative, for any choice of tamed almost-complex structure~$\tilde{J}$. By lemma~\ref{LemmePrincipal}, we only need to prove that for $\delta$ sufficiently small, $[\Sigma]\in\E_{\min}(\widetilde{M}_{\delta})$. But this follows from Lemma~\ref{MinimalExceptionalClass}.
\end{proof}

\begin{example}
Combining Proposition~\ref{ManifoldsWithoutActions} with Theorem~\ref{KedraGeneralise}, one gets examples of rational symplectic $4$-manifolds $(M,\om)$ which do not admit any Hamiltonian action by tori, but whose symplectomorphism groups have infinitely generated rational cohomology. In particular, these topological groups are not homotopy equivalent to homotopy colimits of Lie subgroups, in contrast to $\Symp(S^{2}\times S^{2})$, $\Symp(\PbP)$, and $\Symp(\CP^{2}\nblowup{2})$, see~\cite{AGK,Pi}\eoe
\end{example}

%%%%%%%%%%%%%%%%%%%%%%%%%%%%%%%%%%%%%%%%%%%%%%%%%%%%%%%%%%%%%%%%%%%%%%%%%%%%%%%%
%: REFERENCES
%%%%%%%%%%%%%%%%%%%%%%%%%%%%%%%%%%%%%%%%%%%%%%%%%%%%%%%%%%%%%%%%%%%%%%%%%%%%%%%%

\end{document}